\documentclass[12pt,vatola]{article}

\def\c{\centerline}

\def\re#1{\par\hangindent\parindent\indent\llap{#1\enspace}\ignorespaces}

\def\no{\noindent}

\begin{document}

\c{\large\bf  On Algebraic Multi-Group Spaces}

\vskip 8mm

\c{Linfan Mao}

\vskip 3mm

\c{\scriptsize Institute of Systems Science of Academy of
Mathematics and Systems}

\c{\scriptsize Chinese Academy of Sciences, Beijing 100080,
P.R.China}

\c{\scriptsize E-mail: maolinfan@163.com}

\vskip 8mm

\begin{minipage}{130mm}

\no{\bf Abstract}: A Smarandache multi-space is a union of $n$
spaces $A_1,A_2,\cdots ,A_n$ with some additional conditions
holding. Combining classical of a group with Smarandache
multi-spaces, the conception of a multi-group space is introduced
in this paper, which is a generalization of the classical
algebraic structures, such as the group, filed, body, $\cdots$,
etc.. Similar to groups, some characteristics of a multi-group
space are obtained in this paper.

\vskip 2mm

\no{{\bf Key words:} {\small  multi-space, group, multi-group
space, theorem.}}

\vskip 2mm

\no{{\bf Classification:}  AMS(2000) {\small 05C15, 20H15, 51D99,
51M05}}

\end{minipage}

\vskip 10mm

\no{\bf 1.Introduction}

\vskip 5mm

The notion of multi-spaces is introduced by Smarandache in $[5]$
under his idea of hybrid mathematics: {\it  combining different
fields into a unifying field}($[6]$). Today, this idea is widely
accepted by the world of sciences. For mathematics, definite or
exact solution under a given condition is not the only object for
mathematician. New creation power has emerged and new era for the
mathematics has come now.

A Smarandache multi-space is defined by

\vskip 4mm

\no{\bf Definition $1.1$} \ {\it For any integer $i, 1\leq i\leq
n$ let $A_i$ be a set with ensemble of law $L_i$, and the
intersection of $k$ sets $A_{i_1},A_{i_2},\cdots , A_{i_k}$ of
them constrains the law $I(A_{i_1},A_{i_2},\cdots , A_{i_k})$.
Then the union of $A_i$, $1\leq i\leq n$

$$\widetilde{A} \ = \ \bigcup\limits_{i=1}^n A_i$$

\no is called a multi-space.}\vskip 3mm

The conception of multi-group space is a generalization of the
classical algebraic structures, such as the group, filed, body,
$\cdots$, etc., which is defined as follows.

\vskip 4mm

\no{\bf Definition $1.2$} \ {\it Let
$\widetilde{G}=\bigcup\limits_{i=1}^n G_i$ be a complete
multi-space with a binary operation set
$O(\widetilde{G})=\{\times_i, 1\leq i\leq n\}$. If for any integer
$i, 1\leq i\leq n$, $(G_i;\times_i)$ is a group and for $\forall
x,y,z\in \widetilde{G}$ and any two binary operations ¡°$\times$¡±
and ¡°$\circ$¡±, $\times\not= \circ$, there is one operation, for
example the operation $\times$ satisfying the distribution law to
the operation ¡°$\circ$¡± if their operation results exist , i.e.,

$$x\times (y\circ z) = (x\times y)\circ (x\times z),$$

$$(y\circ z)\times x = (y\times x)\circ (z\times x),$$

\no then $\widetilde{G}$ is called a multi-group space.}

\vskip 3mm

\no{\bf Remark}: \ The following special cases convince us that
multi-group spaces are generalization of group, field and body,
$\cdots$, etc..\vskip 2mm

($i$) If $n=1$, then $\widetilde{G}=(G_1;\times_1)$ is just a
group.

($ii$) If $n=2, G_1=G_2=\widetilde{G}$,  Then $\widetilde{G}$ is a
body. If $(G_1;\times_1)$ and $(G_2;\times_2)$ are commutative
groups, then $\widetilde{G}$ is a field.\vskip 2mm

Notice that in $[7][8]$ various bispaces, such as bigroup,
bisemigroup, biquasigroup, biloop, bigroupoid, biring, bisemiring,
bivector, bisemivector, binear-ring, $\cdots$, etc., consider two
operation on two different sets are introduced.

\vskip 8mm

{\bf $2.$ Characteristics of multi-group spaces}

\vskip 5mm

For a multi-group space $\widetilde{G}$ and a subset
$\widetilde{G_1}\subset\widetilde{G}$, if $\widetilde{G_1}$ is
also a multi-group space under a subset $O(\widetilde{G_1}),
O(\widetilde{G_1})\subset O(\widetilde{G})$, then $\widetilde{G}$
is called a {\it multi-group subspace}, denoted by
$\widetilde{G_1} \preceq \widetilde{G}$. We have the following
criterion for the multi-group subspaces.

\vskip 4mm

\no{\bf Theorem $2.1$} \ {\it For a multi-group space
$\widetilde{G}=\bigcup\limits_{i=1}^nG_i$ with an operation set
$O(\widetilde{G})=\{\times_i|1\leq i\leq n\}$, a subset
$\widetilde{G_1}\subset\widetilde{G}$ is a multi-group subspace if
and only if for any integer $k, 1\leq k\leq n$,
$(\widetilde{G_1}\bigcap G_k;\times_k)$ is a subgroup of
$(G_k;\times_k)$ or $\widetilde{G_1}\bigcap G_k=\emptyset$.}

\vskip 3mm

{\it Proof} \ If $\widetilde{G_1}$ is a multi-group space with the
operation set $O(\widetilde{G_1})=\{\times_{i_j}|1\leq j\leq
s\}\subset O(\widetilde{G})$, then

$$\widetilde{G_1}=\bigcup\limits_{i=1}^n(\widetilde{G_1}\bigcap
G_i) =\bigcup\limits_{j=1}^sG'_{i_j}$$

\no where $G'_{i_j}\preceq G_{i_j}$ and $(G_{i_j};\times_{i_j})$
is a group. Whence, if $\widetilde{G_1}\bigcap G_k\not=\emptyset$,
then there exist an integer $l, k=i_l$ such that
$\widetilde{G_1}\bigcap G_k=G'_{i_l}$, i.e.,
$(\widetilde{G_1}\bigcap G_k;\times_k)$ is a subgroup of
$(G_k;\times_k)$.

Now if for any integer $k$, $(\widetilde{G_1}\bigcap
G_k;\times_k)$ is a subgroup of $(G_k;\times_k)$ or
$\widetilde{G_1}\bigcap G_k=\emptyset$, let $N$ denote the index
set $k$ with $\widetilde{G_1}\bigcap G_k\not=\emptyset$. Then

$$\widetilde{G_1}=\bigcup\limits_{j\in N}(\widetilde{G_1}\bigcap
G_j)$$

\no and $(\widetilde{G_1}\bigcap G_j,\times_j)$ is a group. Since
$\widetilde{G_1}\subset\widetilde{G}$ and
$O(\widetilde{G_1})\subset O(\widetilde{G})$, the associative law
and distribute law are true for the $\widetilde{G_1}$. Therefore,
$\widetilde{G_1}$ is a multi-group subspace of
$\widetilde{G}$.\quad\quad $\natural$

For a finite multi-group subspace, we have the following
criterion.

\vskip 4mm

\no{\bf Theorem $2.2$} \ {\it Let $\widetilde{G}$ be a finite
multi-group space with an operation set
$O(\widetilde{G})=\{\times_i|1\leq i\leq n\}$. A subset
$\widetilde{G_1}$ of $\widetilde{G}$ is a multi-group subspace
under an operation subset $O(\widetilde{G_1})\subset
O(\widetilde{G})$ if and only if for each operation ¡°$\times$¡±
in $O(\widetilde{G_1})$, $(\widetilde{G_1};\times)$ is complete.}

\vskip 3mm

{\it Proof} \ Notice that for a multi-group space $\widetilde{G}$,
its each multi-group subspace $\widetilde{G_1}$ is complete.

Now if $\widetilde{G_1}$ is a complete set under each operation
¡°$\times_i$¡± in $O(\widetilde{G_1})$, we know that
$(\widetilde{G_1}\bigcap G_i; \times_i)$ is a group (see also
$[9]$) or an empty set. Whence, we get that

$$\widetilde{G_1}=\bigcup\limits_{i=1}^n(\widetilde{G_1}\bigcap
G_i).$$

\no Therefore, $\widetilde{G_1}$ is a multi-group subspace of
$\widetilde{G}$ under the operation set $O(\widetilde{G_1})$.
\quad\quad $\natural$

For a multi-group subspace $\widetilde{H}$ of the multi-group
space $\widetilde{G}$, $g\in\widetilde{G}$, define

$$g\widetilde{H}=\{g\times h|h\in\widetilde{H}, \times\in O(\widetilde{H})\}.$$

Then for $\forall x,y\in\widetilde{G}$,

$$x\widetilde{H}\bigcap y\widetilde{H}=\emptyset  \ \ {\rm or} \ \ x\widetilde{H}=y\widetilde{H}.$$

\no In fact, if $x\widetilde{H}\bigcap
y\widetilde{H}\not=\emptyset$, let $z\in x\widetilde{H}\bigcap
y\widetilde{H}$, then there exist elements
$h_1,h_2\in\widetilde{H}$ and operations ¡°$\times_i$¡± and
¡°$\times_j$¡± such that

$$z = x\times_i h_1 = y\times_j h_2.$$

Since $\widetilde{H}$ is a multi-group subspace,
$(\widetilde{H}\bigcap G_i;\times_i)$ is a subgroup. Whence, there
exists an inverse element $h_1^{-1}$ in $(\widetilde{H}\bigcap
G_i;\times_i)$. We have that

$$x\times_i h_1\times_i h_1^{-1} = y\times_j h_2\times_ih_1^{-1}.$$

\no That is,

$$x=y\times_j h_2\times_ih_1^{-1}.$$

\no Whence,

$$x\widetilde{H}\subseteq y\widetilde{H}.$$

Similarly, we can also get that

$$x\widetilde{H}\supseteq y\widetilde{H}.$$

\no Therefore, we get that

$$x\widetilde{H} = y\widetilde{H}.$$

Denote the union of two set $A$ and $B$ by $A\bigoplus B$ if
$A\bigcap B=\emptyset$. Then we get the following result by the
previous proof.

\vskip 4mm

\no{\bf Theorem $2.3$} \ {\it For any multi-group subspace
$\widetilde{H}$ of a multi-group space $\widetilde{G}$, there is a
representation set $T$, $T\subset\widetilde{G}$, such that }

$$\widetilde{G}=\bigoplus\limits_{x\in T}x\widetilde{H}.$$

\vskip 3mm

For the case of finite groups, since there is only one binary
operation ¡°$\times$¡± and $|x\widetilde{H}|=|y\widetilde{H}|$ for
any $x,y\in\widetilde{G}$, We get the following corollary, which
is just Lagrange theorem for finite groups.

\vskip 4mm

\no{\bf Corollary $2.1$}(Lagrange theorem) \ {\it For any finite
group $G$, if $H$ is a subgroup of $G$, then $|H|$ is a divisor of
$|G|$.}

\vskip 3mm

For a multi-group space $\widetilde{G}$ and $g\in\widetilde{G}$,
denote by $\overrightarrow{O(g)}$ all the binary operations
associative with $g$ and by $\widetilde{G}(\times)$ the elements
associative with the binary operation ¡°$\times$¡±. For a
multi-group subspace $\widetilde{H}$ of $\widetilde{G}$,
$\times\in O(\widetilde{H})$ and $\forall
g\in\widetilde{G}(\times)$, if $\forall h\in\widetilde{H}$,

$$g\times h\times g^{-1}\in\widetilde{H},$$

\no then call $\widetilde{H}$ a {\it normal multi-group subspace}
of $\widetilde{G}$, denoted by
$\widetilde{H}\triangleleft\widetilde{G}$. If $\widetilde{H}$ is a
normal multi-group subspace of $\widetilde{G}$, similar to the
normal subgroup of a group, it can be shown that
$g\times\widetilde{H}=\widetilde{H}\times g$, where
$g\in\widetilde{G}(\times)$. We have the following result.

\vskip 4mm

\no{\bf Theorem $2.4$} \ {\it Let
$\widetilde{G}=\bigcup\limits_{i=1}^nG_i$ be a multi-group space
with an operation set $O(\widetilde{G})=\{\times_i|1\leq i\leq
n\}$. Then a multi-group subspace $\widetilde{H}$ of
$\widetilde{G}$ is normal if and only if for any integer $i, 1\leq
i\leq n$, $(\widetilde{H}\bigcap G_i;\times_i)$ is a normal
subgroup of $(G_i;\times_i)$ or $\widetilde{H}\bigcap
G_i=\emptyset$.}

\vskip 3mm

{\it Proof} \ We have known that

$$\widetilde{H} = \bigcup\limits_{i=1}^n(\widetilde{H}\bigcap G_i).$$

If for any integer $i, 1\leq i\leq n$, $(\widetilde{H}\bigcap
G_i;\times_i)$ is a normal subgroup of $(G_i;\times_i)$, then we
know that for $\forall g\in G_i, 1\leq i\leq n$,

$$g\times_i(\widetilde{H}\bigcap G_i)\times_ig^{-1} = \widetilde{H}\bigcap G_i.$$

\no Whence, for $\forall \circ\in O(\widetilde{H})$ and $\forall
g\in\overrightarrow{\widetilde{G}(\circ )}$,

$$g\circ\widetilde{H}\circ g^{-1} = \widetilde{H}.$$

\no That is, $\widetilde{H}$ is a normal multi-group subspace of
$\widetilde{G}$.

Now if $\widetilde{H}$ is a normal multi-group subspace of
$\widetilde{G}$, then by definition, we know that for $\forall
\circ\in O(\widetilde{H})$ and $\forall g\in\widetilde{G}(\circ
)$,

$$g\circ\widetilde{H}\circ g^{-1} = \widetilde{H}.$$

Not loss of generality, we assume that $\circ = \times_k$, then we
get that

$$g\times_k(\widetilde{H}\bigcap G_k)\times_k g^{-1}= \widetilde{H}\bigcap G_k.$$

\no Therefore, $(\widetilde{H}\bigcap G_k; \times_k)$ is a normal
subgroup of $(G_k,\times_k)$. For operation ¡°$\circ$¡± is chosen
arbitrarily, we know that for any integer $i$, $1\leq i\leq n$,
$(\widetilde{H}\bigcap G_i;\times_i)$ is a normal subgroup of
$(G_i;\times_i)$ or an empty set.\quad\quad $\natural$

For a multi-group space $\widetilde{G}$ with an operation set
$O(\widetilde{G})=\{\times_i| \ 1\leq i\leq n\}$, an order of
operations in $O(\widetilde{G})$ is said an {\it oriented
operation sequence}, denoted by
$\overrightarrow{O}(\widetilde{G})$. For example, if
$O(\widetilde{G})=\{\times_1,\times_2\times_3\}$, then
$\times_1\succ\times_2\succ\times_3$ is an oriented operation
sequence and $\times_2\succ\times_1\succ\times_3$ is also an
oriented operation sequence.

For an oriented operation sequence
$\overrightarrow{O}(\widetilde{G})$, we construct a series of
normal multi-group subspaces

$$\widetilde{G}\triangleright\widetilde{G}_1\triangleright\widetilde{G}_2
\triangleright\cdots\triangleright\widetilde{G}_m=\{1_{\times_n}\}$$

\no by the following programming.

\vskip 3mm

\no{\bf STEP $1$}: {\it Construct a series

$$\widetilde{G}\triangleright\widetilde{G}_{11}\triangleright\widetilde{G}_{12}
\triangleright\cdots\triangleright\widetilde{G}_{1l_1}$$

\no under the operation ¡°$\times_1$¡±.}

\vskip 2mm

\no{\bf STEP $2$}: {\it If a series

$$\widetilde{G}_{(k-1)l_1}\triangleright\widetilde{G}_{k1}\triangleright\widetilde{G}_{k2}
\triangleright\cdots\triangleright\widetilde{G}_{kl_k}$$

\no has be constructed under the operation ¡°$\times_k$¡± and
$\widetilde{G}_{kl_k}\not=\{1_{\times_n}\}$, then construct a
series

$$\widetilde{G}_{kl_1}\triangleright\widetilde{G}_{(k+1)1}\triangleright\widetilde{G}_{(k+1)2}
\triangleright\cdots\triangleright\widetilde{G}_{(k+1)l_{k+1}}$$

\no under the operation ¡°$\times_{k+1}$¡±.

This programming is terminated until the series

$$\widetilde{G}_{(n-1)l_1}\triangleright\widetilde{G}_{n1}\triangleright\widetilde{G}_{n2}
\triangleright\cdots\triangleright\widetilde{G}_{nl_n}=\{1_{\times_n}\}$$

\no has be constructed under the operation ¡°$\times_n$¡±.}

The number $m$ is called the length of the series of normal
multi-group subspaces. For a series

$$\widetilde{G}\triangleright\widetilde{G}_1\triangleright\widetilde{G}_2
\triangleright\cdots\triangleright\widetilde{G}_n=\{1_{\times_n}\}$$

\no of normal multi-group subspaces, if for any integer $k, s,
1\leq k\leq n, 1\leq s\leq l_k$, there exists a normal multi-group
subspace $\widetilde{H}$ such that

$$\widetilde{G}_{ks}\triangleright\widetilde{H}\triangleright\widetilde{G}_{k(s+1)},$$

\no then $\widetilde{H}=\widetilde{G}_{ks}$ or
$\widetilde{H}=\widetilde{G}_{k(s+1)}$, we call this series is
{\it maximal}. For a maximal series of finite normal multi-group
subspaces, we have the following result.

\vskip 4mm

\no{\bf Theorem $2.5$} \ {\it For a finite multi-group space
$\widetilde{G}=\bigcup\limits_{i=1}^n G_i$ and an oriented
operation sequence $\overrightarrow{O}(\widetilde{G})$, the length
of maximal series of normal multi-group subspaces is a constant,
only dependent on $\widetilde{G}$ itself.}

\vskip 3mm

{\it Proof} \ The proof is by induction on the integer $n$.

For $n=1$, the maximal series of normal multi-group subspaces is
just a composition series of a finite group. By Jordan-H\"{o}lder
theorem (see $[1]$ or $[3]$), we know the length of a composition
series is a constant, only dependent on $\widetilde{G}$. Whence,
the assertion is true in the case of $n=1$.

Assume the assertion is true for cases of $n\leq k$. We prove it
is true in the case of $n=k+1$. Not loss of generality, assume the
order of binary operations in $\overrightarrow{O}(\widetilde{G})$
being $\times_1\succ\times_2\succ\cdots\succ\times_n$ and the
composition series of the group $(G_1,\times_1)$ being

$$G_1\triangleright G_2\triangleright\cdots\triangleright G_s=\{1_{\times_1}\}.$$

By Jordan-H\"{o}lder theorem, we know the length of this
composition series is a constant, dependent only on
$(G_1;\times_1)$. According to Theorem $3.6$, we know a maximal
series of normal multi-group subspace of $\widetilde{G}$ gotten by
the STEP $1$ under the operation ¡°$\times_1$¡± is

$$\widetilde{G}\triangleright\widetilde{G}\setminus (G_1\setminus G_2)\triangleright
\widetilde{G}\setminus (G_1\setminus
G_3)\triangleright\cdots\triangleright\widetilde{G}\setminus
(G_1\setminus \{1_{\times_1}\}).$$

Notice that $\widetilde{G}\setminus
(G_1\setminus\{1_{\times_1}\})$ is still a multi-group space with
less or equal to $k$ operations. By the induction assumption, we
know the length of its maximal series of normal multi-group
subspaces is only dependent on $\widetilde{G}\setminus
(G_1\setminus\{1_{\times_1}\})$, is a constant. Therefore, the
length of a maximal series of normal multi-group subspaces is also
a constant, only dependent on $\widetilde{G}$.

Applying the induction principle, we know that the length of a
maximal series of normal multi-group subspaces of $\widetilde{G}$
is a constant under an oriented operations
$\overrightarrow{O}(\widetilde{G})$, only dependent on
$\widetilde{G}$ itself.\quad\quad $\natural$\vskip 2mm

As a special case, we get the following corollary.

\vskip 4mm

\no{\bf Corollary $2.2$}(Jordan-H\"{o}lder theorem) \ {\it For a
finite group $G$, the length of the composition series is a
constant, only dependent on $G$.}

\vskip 8mm

{$3.$ Open Problems on Multi-group Spaces}

\vskip 5mm

\no{\bf Problem $3.1$} \ {\it Establish a decomposition theory for
multi-group spaces.}

\vskip 2mm

In group theory, we know the following decomposition
results([$1$][$3$]) for a group.

\vskip 2mm

{\it Let $G$ be a finite $\Omega$-group. Then $G$ can be uniquely
decomposed as a direct product of finite non-decomposition
$\Omega$-subgroups.}

\vskip 2mm

{\it Each finite Abelian group is a direct product of its Sylow
$p$-subgroups.}\vskip 2mm

Then Problem $3.1$ can be restated as follows.

\vskip 2mm

{\it Whether can we establish a decomposition theory for
multi-group spaces similar to above two results in group theory,
especially, for finite multi-group spaces?}

\vskip 3mm

\no{\bf Problem $3.2$} \ {\it Define the conception of simple
multi-group spaces for multi-group spaces. For finite multi-group
spaces, whether can we find all simple multi-group spaces?}\vskip
2mm

For finite groups, we know that there are four simple group
classes ([$9$]):

\vskip 2mm

{\bf Class $1$}: the cyclic groups of prime order; \vskip 2mm

{\bf Class $2$}:the alternating groups $A_n, n\geq 5$;

\vskip 2mm

{\bf Class $3$}: the 16 groups of Lie types;

\vskip 2mm

{\bf Class $4$}: the 26 sporadic simple groups.

\vskip 3mm

\no{\bf Problem $2.3$} \ {\it Determine the structure properties
of a multi-group space generated by finite elements.} \vskip 2mm

For a subset $A$ of a multi-group space $\widetilde{G}$, define
its spanning set by

$$\left<A\right>=\{a\circ b | a,b\in A \ {\rm and} \ \circ\in
O(\widetilde{G})\}.$$

\no If there exists a subset $A\subset\widetilde{G}$ such that
$\widetilde{G}=\left<A\right>$, then call $\widetilde{G}$ is
generated by $A$. Call $\widetilde{G}$ is {\it finitely generated}
if there exist a finite set $A$ such that
$\widetilde{G}=\left<A\right>$. Then Problem $2.3$ can be restated
by \vskip 2mm

{\it Can we establish a finite generated multi-group theory
similar to the finite generated group theory?}

\vskip 8mm

{\bf References}

\vskip 5mm

\re{[1]}G.Birkhoff and S.Mac Lane, {\it A Survey of Modern
Algebra}, Macmillan Publishing Co., Inc, 1977.

\re{[2]}Daniel Deleanu, {\it A Dictionary of Smarandache
Mathematics}, Buxton University Press, London \& New York,2004.

\re{[3]}Lingzhao Nie and Shishun Ding, {\it Introduction to
Algebra}, Higher Education Publishing Press, 1994.

\re{[4]}L.F.Mao, {\it Automorphism Groups of Maps, Surfaces and
Smarandache Geometries}, American Research Press, 2005.

\re{[5]}F.Smarandache, Mixed noneuclidean geometries, {\it eprint
arXiv: math/0010119}, 10/2000.

\re{[6]}F.Smarandache, {\it A Unifying Field in Logics.
Neutrosopy: Neturosophic Probability, Set, and Logic}, American
research Press, Rehoboth, 1999.

\re{7]]}W.B.Vasantha Kandasamy, {\it Bialgebraic structures and
Smarandache bialgebraic structures}, American Research Press,
2003.

\re{[8]}W.B.Vasantha Kandasamy and F.Smarandache, {\it Basic
Neutrosophic Algebraic Structures and Their Applications to Fuzzy
and Neutrosophic Models}, HEXIS, Church Rock, 2004.

\re{[9]}Mingyao Xu, {\it Introduction to Group Theory}(I)(II),
Science Publish Press, Beijing ,1999.

\end{document}